\documentclass{IEEEtran}
\usepackage[utf8]{inputenc}

\IEEEoverridecommandlockouts

\usepackage{cite}
\usepackage{algorithmic}
\usepackage{textcomp}

\usepackage{graphicx}
\usepackage{amsmath}
\usepackage{amssymb}
\usepackage{amsfonts}
\usepackage{mathrsfs}
\usepackage{dsfont}
\makeatletter
\def\amsbb{\use@mathgroup \M@U \symAMSb}
\makeatother
\usepackage{color}

\usepackage[usenames,dvipsnames]{xcolor}
\usepackage{epstopdf}
\usepackage{siunitx}

\newtheorem{theorem}{Theorem}
\newtheorem{remark}{Remark}
\newtheorem{definition}{Definition}
\newtheorem{lemma}{Lemma}
\newtheorem{proposition}{Proposition}
\newtheorem{corollary}{Corollary}
\newtheorem{assumption}{Assumption}
\newtheorem{design condition}{Design condition}
\newcommand{\vect}[1]{\boldsymbol{#1}} 
\newcommand{\norm}[1]{\left\lVert#1\right\rVert}
\usepackage{enumerate}
\usepackage{balance}

\usepackage{siunitx}

\usepackage{tabularx,booktabs,caption}

\usepackage{mathtools}

\DeclarePairedDelimiter\floor{\lfloor}{\rfloor}

\allowdisplaybreaks

\newcommand{\iclf}[1]{{\color{black} #1}}

\newcommand{\il}[1]{{\color{black} #1}}
\newcommand{\an}[1]{{\color{black} #1}}
\newcommand{\ill}[1]{{\color{black} #1}}
\newcommand{\ank}[1]{{\color{black} #1}}
\newcommand{\illl}[1]{{\color{black} #1}}
\newcommand{\ka}[1]{{\color{black} #1}}
\newcommand{\icll}[1]{{\color{black} #1}}
\newcommand{\kas}[1]{{\color{black} #1}}
\newcommand{\ic}[1]{{\color{black} #1}}
\newcommand{\icc}[1]{{\color{black} #1}}
\newcommand{\iccc}[1]{{\color{black} #1}}
\newcommand{\ankk}[1]{{\color{black} #1}}
\newcommand{\ilc}[1]{{\color{black} #1}}
\newcommand{\ilcc}[1]{{\color{black} #1}}
\newcommand{\ant}[1]{{\color{black} #1}}
\newcommand{\antk}[1]{{\color{black} #1}}

\newcommand{\antkk}[1]{{\color{black} #1}}

\title{Frequency \ant{R}egulation \il{with} \\\ant{T}hermostatically \ant{C}ontrolled \iclf{\ant{L}oads}:\\ \il{\ant{A}ggregation of \ant{D}ynamics and \ant{S}ynchronization}}

\author{Andreas Kasis\thanks{\antkk{Andreas Kasis is with the KIOS Research and Innovation Center of Excellence and the Department of Electrical and Computer Engineering, University of Cyprus, Cyprus; e-mail: kasis.andreas@ucy.ac.cy}}
 and Ioannis Lestas\thanks{Ioannis Lestas is with the Department of Engineering, University of Cambridge,  Cambridge,  United Kingdom; e-mail: icl20@cam.ac.uk}
 \thanks{{This work was supported by ERC starting grant No. 679774 and by the European Union’s Horizon 2020 research and innovation program under grant agreement No. 891101 (SmarTher Grid).}}
\thanks{\ka{A preliminary version of this work has appeared in \cite{kasis2019_cdc}. This manuscript \kas{extends the analysis to a \ic{broad} class of linear generation dynamics and includes} the analytic proofs of the main results, additional discussion and simulations that demonstrate the impact of the proposed analysis.}}
}

\begin{document}

\maketitle

\begin{abstract}
Thermostatically controlled loads (TCLs) can provide ancillary services to the power network by aiding existing frequency control mechanisms.
{TCLs} are, {however}, characterized by {an} intrinsic limit cycle behavior which raises the risk that these could synchronize
when coupled with {the frequency dynamics of the power grid}, i.e. simultaneously switch, inducing persistent and possibly catastrophic power oscillations.
{To address this problem, schemes
{with a randomized response time in their control policy}
 have been proposed in the literature.}
{However,  {such schemes}  introduce delays in the response of TCLs {to frequency feedback} that may limit their ability to {provide {fast} support at urgencies.}}
In this paper, we present a deterministic control mechanism for TCLs such that those switch when prescribed frequency thresholds are exceeded in order to provide ancillary services to the power network. 
For the considered scheme,
{we provide analytic conditions which ensure that synchronization is avoided.}
{In particular, we show that} as the number of loads tends to infinity,
there exist arbitrarily long time intervals where the frequency deviations are arbitrarily {small. 
Our {analytical}} results are verified with simulations on the Northeast Power Coordinating Council (NPCC) 140-bus system, which demonstrate that the proposed scheme offers {improved frequency response compared with existing implementations.}
\end{abstract}

\section{Introduction}
\textbf{Motivation and literature review:}
A significant growth in {the} penetration of renewable sources of generation in power networks is expected {over} the following years \cite{lund2006large, ipakchi2009grid},  driven by environmental concerns.
This will result in  increasingly intermittent generation, endangering power quality and potentially  the stability of the power network.
Controllable loads are considered to be a way to counterbalance intermittent generation, due to their ability to provide a fast response at urgencies by accordingly adapting their demand.
The use of loads as ancillary services, in conjunction with a large penetration of renewable sources of generation will significantly increase the number of active devices in the network making its electromechanical {response} difficult to predict and encouraging the analytical study of its behavior.
Along these lines,   various research studies in {recent years have} considered controllable demand as a means {of providing support to} primary~\cite{molina2011decentralized, trip2016optimal, kasis2016primary, devane2016primary, kasis2017primary}, and secondary~\cite{mallada2017optimal},  \cite{kasis2017stability, kasis2019secondary}, frequency control mechanisms, {where the objective is} to ensure that generation and demand are balanced and that {the} frequency converges to its nominal value (50Hz or 60Hz) respectively.

Thermostatically controlled loads (TCLs) comprise a significant portion of the total demand. A recent survey   in the EU \cite{zimmermann2012household}  showed that TCLs exceeded $80\%$ and $40\%$ of the total consumption in households with and without electric heating respectively.
{TCLs have an intrinsic limit cycle behavior whereby they need to periodically turn on and off in order to maintain the temperature within a prescribed range. This significantly complicates their use for frequency control, in comparison with loads that are not thermostatically controlled \cite{kasis2019secondary, liu2016non, kasis2021primary}.  In particular,   the coupling of the individual limit cycles in TCLs with the grid frequency, could lead to a synchronization of these limit cycles thus resulting to highly undesirable oscillations in the aggregate load profile. Therefore dedicated analysis tools and studies are needed for the efficient integration of TCLs to the grid such that they provide support to frequency regulation.}

The use of TCLs {for} frequency control has been considered in  \cite{short2007stabilization}, where the authors suggested  temperature thresholds in TCLs to be  linearly dependent on frequency and demonstrated with simulations that this resulted in improved performance.
 However, it was demonstrated in \cite{angeli2012stochastic} that such control {schemes} could potentially result to {load synchronization.}
 As a remedy to this problem, the authors proposed a {randomized} control scheme which ensured that {TCLs would} not synchronize.
{Various other studies considered similar problems by proposing schemes with randomization in the control policy.
{In \cite{aunedi2013economic},  safety constraints in the operation of {TCLs} are additionally included, and \cite{tindemans2015decentralized, totu2017demand} incorporate stochastic switching in the TCL operation so as to achieve a prescribed power profile.
{However, schemes with a randomized delay in their control policies may limit the ability of TCLs to respond to unforeseen frequency fluctuations and provide ancillary support at fast timescales.}}}
The latter, motivates the study of {alternative schemes} for the control of thermostatic loads, such that a faster response {can be achieved} at {urgencies, while at the same time avoiding load synchronization.}

\textbf{Contribution:} This study considers a deterministic approach for the control of thermostatic loads, such that ancillary services {with a fast response} are provided at urgencies.
Our main analytic results concern the case  {where} the number of loads tends to infinity, a condition  justified by the large number of thermostatic appliances in power {networks.}

{More precisely, we propose a} control scheme for TCLs, such that loads switch when certain frequency thresholds are exceeded in order to support existing secondary frequency control schemes.
{For the considered scheme, we provide design conditions for the frequency thresholds that bound the coupling between {the frequency and the} load dynamics {so as to avoid load synchronization}.
{In particular, one of the main results is to analytically show} {that} when the number of loads tends to infinity, {the} frequency deviations will be arbitrarily small for arbitrarily long time intervals.

The proposed scheme also ensures that load temperatures will not exceed their respective bounds, and hence that user comfort levels will not be affected.
Furthermore, the fact that loads switch instantly at urgencies,  leads to a fast response whereby randomized delays, often used in the literature to avoid synchronization, are avoided.}

Our {analytical} results are verified with numerical simulations on the NPCC 140-bus network, where it is demonstrated that  the proposed scheme offers reduced frequency overshoots  in comparison with existing implementations.

\textbf{Paper structure:} In Section \ref{Notation} we present  some basic notation used in the paper and in Section \ref{sec:Network_model} the considered power system.
 In Section \ref{Sec:Conventional_thermostatic} we consider a conventional model for TCLs and study its properties in terms of {the} aggregate mean and variance.
In Section \ref{Sec:Novel_thermostatic}, we present our proposed scheme for frequency control using TCLs and state our main results regarding the performance of the power system.
Numerical investigations of the results on the NPCC 140-bus system are provided  in Section \ref{Simulation} and conclusions are drawn in Section \ref{Conclusion}. The proofs of the main results are provided in  the appendix.

\section{Notation}\label{Notation}

Real, {natural and complex} numbers are denoted by $\mathbb{R}$, $\mathbb{N}$ and {$\mathbb{C}$ respectively}, and the set of n-dimensional vectors with real entries is denoted by $\mathbb{R}^n$. Furthermore, we define the sets of integers and {strictly} positive rational and {strictly} positive real numbers
  by $\mathbb{Z}, \mathbb{Q}_+$ and $\mathbb{R}_+$ respectively.
 The set of natural numbers including zero {is} denoted by $\mathbb{N}_0$.
 The cardinality of a set $S$ is denoted by $|S|$.
 For $a \in \mathbb{R},b \in \mathbb{R} \setminus \{0\}$,
\ilc{$a$~modulo $b$}
is denoted by
$[a]^+_b$ and defined as
$
[a]_{b}^+ = a - b \floor{\frac{a}{b}}
$, where for $x \in \mathbb{R}$, $\floor{x} = \sup \{m \in \mathbb{Z}: m \leq x\}$.
 The average of a real {valued} time signal $x(t)$ with respect to time is {defined as $\mathbb{E}(x(t)) = \lim_{\tau \rightarrow \infty} \frac{1}{\tau} \int_{0}^{\tau} x (t)dt$ and its variance as $\mathbb{V}(x(t)) = \mathbb{E}((x(t))^2) - [\mathbb{E}(x(t))]^2$.}
{{For $c\in\mathbb{C}$ we denote its magnitude by $|c|$.}}
{{The 1-norm of} a linear system with {transfer function $G(s)$}
is given by $\int_{0}^{\infty} |g(t)| dt$, where $g(t)$ is the inverse Laplace transformation of $G(s)$.}
{We use} $\vect{0}_n$ to denote  the $n \times 1$ vector with all elements equal to $0$.
{We also say that a matrix $A\in\mathbb{R}^{n\times n}$ is {Hurwitz} if all its eigenvalues have strictly negative real part.}
{Finally, a sequence $\{s_1, s_2, s_3, ...\}$ of real numbers is said to be uniformly distributed on an interval $[a, b]$ if for any subinterval $[c, d]$ of $[a, b]$ we have
$\lim_{n \rightarrow \infty} \frac{|\{s_1, s_2, s_3, ...\} \cap [c, d]|}{n} = \frac{d-c}{b-a}.$}

\section{Power system model}\label{sec:Network_model}

We use the swing {equation} to describe the rate of change of {the} frequency of the power system  (e.g. \cite{Bergen_Vittal}).
In particular, we consider the following assumptions on our studied model: \newline
1) Bus voltage magnitudes satisfy $|V| = 1$ p.u. for all buses. \newline
2) Lines are lossless and characterized by their susceptances. \newline
3) Reactive power flows do not affect bus voltage phase angles and frequencies.\newline
4) Frequencies between buses are synchronized.

The first three conditions have been widely used in the literature for the study of frequency control schemes in power networks \cite{trip2016optimal, kasis2019secondary}.
The fourth assumption is justified from the relatively small deviations between bus frequencies, which allows the study of power system characteristics using a single frequency (see also \cite{angeli2012stochastic, anderson1990low, shi2018analytical}).
{The latter follows from the fact that the {dynamic} behavior of TCLs is much slower than the frequency dynamics between buses,  which justifies the assumption of 
{small deviations among bus frequencies.}
Please note that a full complexity
{power network model, which includes {multiple buses}, voltage dynamics, line resistances and reactive power flows,} is considered in the simulations presented {in} Section \ref{Simulation}, which verify the main results of the paper.
The above motivate} the following system dynamics,
\begin{equation}
 M \dot{\omega} = - p^L + p^M - D\omega - \sum_{j \in N} d^c_j.
  \label{sys1}
 \end{equation}

In system~\eqref{sys1} the time-dependent variables  $p^M$, $d^c_j$ and $\omega$ represent, respectively,
the {aggregate} mechanical power injection, the $j$th thermostatic load and
the deviation from the nominal value\footnote{{The nominal value is 50Hz or 60Hz.}}
 of {the} frequency. Furthermore, we let $N~:=~\{1,2\dots,|N|\}$ be the set of TCLs.
The constants $M > 0$ and $D > 0$ denote the generator inertia  and damping coefficient respectively.
{Finally, the {aggregate} uncontrollable demand is denoted by $p^L$.}

\subsection{Generation Dynamics}

{We consider a broad class of linear generation dynamics of the form}
\begin{equation} \label{sys2}
p^M = \hat{C} \hat{x} + \hat{D} \omega, \quad \dot{\hat{x}} = \hat{A} \hat{x} + \hat{B} \omega,
\end{equation}
{{with input $\omega$, output $p^M$, state $\hat{x}$ that takes values in $\mathbb{R}^{n}$ and corresponding matrices $\hat{A} \in \mathbb{R}^{n \times n}, \hat{B} \in \mathbb{R}^{n}, \hat{C} \in \mathbb{R}^{1 \times n}$ and $\hat{D} \in \mathbb{R}$.}
Note that linear  {systems} are widely used in the literature to model generation {dynamics} (see e.g.  \cite[Section 11.1]{Bergen_Vittal}, \cite[Section 11.1.7]{machowski2011power}). Such models are particularly relevant when small disturbances are considered.

The system \eqref{sys1}, \eqref{sys2} can be represented {in the} form
\begin{equation} \label{sys1_2}
\begin{bmatrix}
\dot{\omega} \\
\dot{\hat{x}}
\end{bmatrix}
= A\begin{bmatrix}
{\omega} \\
{\hat{x}}
\end{bmatrix}  + B[p^L +   \sum_{j \in \mathcal{N}} d^c_j],
\end{equation}
where
 $A = \begin{bmatrix}
(\hat{D} -D)/M & \hat{C}/M \\
\hat{B} & \hat{A}
\end{bmatrix}$ and
 $B = \begin{bmatrix}
-1/M \\
\vect{0}_n
\end{bmatrix}$. {We also denote $\hat{u}=p^L +   \sum_{j \in \mathcal{N}} d^c_j$.}}
The following assumption is made for \eqref{sys1_2}.
\begin{assumption}\label{assum_cont_system}
For system \eqref{sys1_2} the following hold
\begin{enumerate}[(i)]
\item $A$ is {Hurwitz},
\item All equilibria of \eqref{sys1_2} {with constant $\hat{u}$} satisfy $\omega^* = 0$.
\end{enumerate}
\end{assumption}
Assumption \ref{assum_cont_system}(i) ensures that \eqref{sys1_2} is an asymptotically stable system. The latter is in line with current implementations where generation dynamics are designed {such that the {power} system is {stable.}}
Assumption \ref{assum_cont_system}(ii) is {associated with the fact that secondary frequency control is implemented,} where the objective is {to recover the frequency to} its nominal value at steady state.

\section{Thermostatically controlled loads}\label{Sec:Conventional_thermostatic}

\begin{figure}[t!]
\centering
\includegraphics[trim = 3mm 0mm 0mm 0mm, scale = 0.6,clip=true]{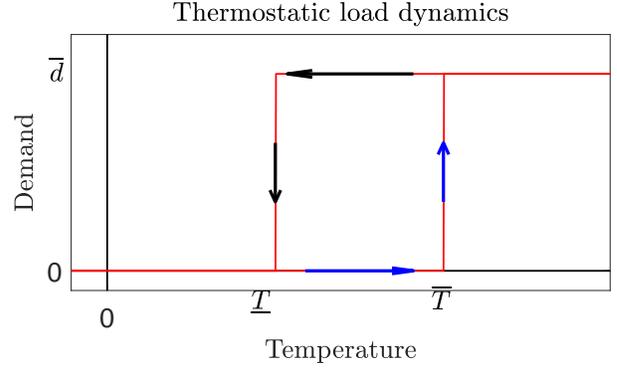}
\vspace{-1mm}
\caption{TCL scheme described by \eqref{sys_hysteresis}.}
\label{Figure_thermo_conventional}
\end{figure}

In this section we consider {a conventional} model for {cooling TCLs (e.g. refrigerators, air conditioner units)}  and study its properties.
{Note that the extension to heating TCLs, such as space heaters, is trivial and thus omitted.}
The analysis below enables to deduce important properties of TCL behavior, which are used to obtain the main results of this paper.
TCL dynamics are {commonly  described by (e.g \cite{angeli2012stochastic, stadler2009modelling})}
\begin{equation}\label{sys_hysteresis}
d^c_j = \overline{d}_j \sigma_{j}, \quad
\sigma_{j}(t^+) = \begin{cases}
1, \qquad T_j \geq  \overline{T}_j,\\[1mm]
0,  \qquad  T_j \leq  \underline{T}_j,\\[1mm]
\sigma_j(t), \hspace{1.8mm}    \underline{T}_j \leq T_j  \leq  \overline{T}_j,
\end{cases} \hspace{-2mm}
\end{equation}
where  $j\in N$ and  $t^+ = \lim_{\epsilon \rightarrow 0} (t + \epsilon)$. In \eqref{sys_hysteresis}, the time-dependent variables $d^c_j$, and ${\sigma_j \in \{0,1\}}$ denote the demand and switch state of the $j$th load respectively. {The time dependent variable $T_j$ denotes the temperature of the $j$th load}.
The constants $\overline{d}_j, \underline{T}_j$ and $\overline{T}_j$ denote the load magnitude and lower and upper temperature thresholds for load $j$ respectively and satisfy
 $\overline{d}_j \in \mathbb{R}_+$ and  $\overline{T}_j > \underline{T}_j >0, j \in N$.
 The hysteresis scheme in \eqref{sys_hysteresis} is depicted in Figure \ref{Figure_thermo_conventional}.

Furthermore, the temperature dynamics satisfy
\begin{equation}\label{sys_temperature}
\dot{T}_j = -k_j (T_j - \hat{T}_j + \lambda_j d^c_j), j \in N,
\end{equation}
where {constants} {$k_j, \lambda_j > 0$ denote the} thermal insulation coefficient and coefficient of performance  of load $j$ respectively. Furthermore,  $\hat{T}_j$ denotes the ambient temperature of load $j$ that is assumed {to be} constant.  Moreover, it is assumed that $\hat{T}_j - \lambda_j \overline{d}_j < \underline{T}_j$ and $\hat{T}_j > \overline{T}_j, j \in N$, such that\footnote{Note that the conditions $\hat{T}_j - \lambda_j \overline{d}_j < \underline{T}_j$ and $\hat{T}_j > \overline{T}_j$ correspond to cooling devices, such as air-conditioning units and refrigerators. These inequalities should be appropriately adapted for heating units, such as space heaters. {This} extension in the analysis is trivial and {is} hence omitted.} \eqref{sys_hysteresis}, \eqref{sys_temperature}, has no equilibria, as is the case in practice.

\subsection{Periods and duty cycles of TCLs}

The period $\pi_j$ of thermal load $j$, described by \eqref{sys_hysteresis}, \eqref{sys_temperature}, is defined as the time required for load $j$  to switch twice, i.e. the time between two consecutive switches to the ON (or equivalently OFF) state.
  In the following definition, we let $t_{j,i}$ be the time where the $i$th switch of load $j$, described by \eqref{sys_hysteresis},~\eqref{sys_temperature}, occurs.

\begin{definition}\label{period_dfn}
The  period of load $j$ is defined as $\pi_j = t_{j,i+2} - t_{j,i}$, for any $i \geq 1$.
\end{definition}

It should be clear that for any $j \in N$, it holds that $t_{j,i+2} - t_{j,i} = t_{j,k+2} - t_{j,k}$, for all $i,k \in \mathbb{N}$.
Note that, as follows from  \eqref{sys_hysteresis}, \eqref{sys_temperature}, the time lengths that load $j$ remains switched ON and OFF  within each period are respectively given by
\begin{subequations}\label{periods_ON_OFF}
\begin{gather}
\pi_j^{ON} = \frac{1}{k_j} \ln(\frac{\overline{T}_j + \lambda_j \overline{d}_j - \hat{T}_j}{\underline{T}_j  + \lambda_j \overline{d}_j - \hat{T}_j} ), j \in N, \\
\pi_j^{OFF} = \frac{1}{k_j} \ln(\frac{\hat{T}_j - \underline{T}_j}{\hat{T}_j - \overline{T}_j} ), j \in N,
\end{gather}
\end{subequations}
and that it  trivially follows that $\pi_j = \pi_j^{ON} + \pi_j^{OFF}$.
Furthermore, the duty cycle of each load is given by $\alpha_j = \frac{\pi_j^{ON}}{\pi_j}$, i.e. the ratio of time the load is ON {within} each period.
 Moreover, we define the period ratio between loads $i$ and $j$ as $\rho_{ij} = \frac{\pi_i}{\pi_j}$.
We shall use $d^{c,*}_j = \alpha_j \overline{d}_j$ to denote the average value of $d^{c}_j$ when its dynamics are described by \eqref{sys_hysteresis},~\eqref{sys_temperature}.
{In addition},  we~{{let}
\begin{equation}
d^s = \sum_{j \in N} d^c_j, \quad \Gamma = \sum_{j \in N} \overline{d}_j,
\end{equation}
be the aggregate sum and aggregate magnitude of TCLs,  where $\Gamma \in \mathbb{R}_+$.}
Finally, we define $E = \{ (i,j) : i,j \in N, i \neq j \}$ as the set of all load pairs.

\subsection{Variance analysis}

In this section we consider the {aggregation} of TCLs {and analyze its} mean and variance. In particular, we study how {the latter is} influenced when the number of loads tends to infinity, assuming {a constant} aggregate sum.

An important assumption in the following analysis is that period ratios lie in the set $\mathbb{R}_+ / \mathbb{Q}_+$.
This is stated below.

\begin{assumption}\label{assumption_period_ratio}
All loads $(i,j) \in E$ described by \eqref{sys_hysteresis},~\eqref{sys_temperature}, satisfy $\rho_{ij} \in \mathbb{R}_+ / \mathbb{Q}_+$.
\end{assumption}

Assumption \ref{assumption_period_ratio} is {a technical} condition that enables to deduce Theorem \ref{thm_variance} below which shows that {when the number of TCLs tends to infinity, then the variance of their aggregation} is zero for any initial condition. {In particular, when Assumption \ref{assumption_period_ratio} holds, {then} $d^s$  is an aperiodic signal that exhibits variability in the time instances the individual loads switch on and off{,} thus leading to Theorem 1.}
Assumption \ref{assumption_period_ratio} excludes cases {where} 
two loads have identical periods, which makes the  aggregation of any two loads periodic. 
The latter is true for all cases where $\rho_{ij} \in \mathbb{Q}_+$, which are hence excluded.
Note that  $\mathbb{Q}_+$ is {a set} of measure zero and hence the condition $\rho_{ij} \in \mathbb{R}_+ / \mathbb{Q}_+$ is unlikely to be violated in practice.

The following theorem states that the variance of the aggregation of TCLs tends to zero as their number tends to infinity. Its proof can be found in the appendix.

\begin{theorem}\label{thm_variance}
Consider thermostatic loads  described by \eqref{sys_hysteresis},~\eqref{sys_temperature}, with $\overline{d}_j = \frac{\Gamma}{|N|}$ and let Assumption \ref{assumption_period_ratio} hold. Then,
{$\mathbb{V}(d^s) {<} \frac{\Gamma^2}{|N|}$ and hence}
 $\lim_{|N| \rightarrow \infty}\mathbb{V}(d^s) = 0$.
\end{theorem}
{Theorem \ref{thm_variance} demonstrates}
that as the number of  loads described by \eqref{sys_hysteresis}, \eqref{sys_temperature},  becomes large, then an almost flat {aggregate {demand}} should be expected, a desired feature {to avoid large oscillations in the frequency response.} 
Note that Theorem \ref{thm_variance}, as well as many of the results that follow,  are stated for the case where $\overline{d}_j = \frac{\Gamma}{|N|}, j \in N,$ which suggests a constant aggregate sum $\Gamma$ and loads of identical magnitude.
{The assumption that all load magnitudes are identical is made  
for simplicity and could potentially be relaxed, as part of future work.}

\begin{remark}\label{remark_stochastic}
A result analogous to Theorem \ref{thm_variance} could be obtained by adopting a stochastic description for TCLs, {where {these} are modeled as independent random processes}.
Theorem \ref{thm_variance} is stated based on the presented deterministic setting, described by \eqref{sys_hysteresis}--\eqref{sys_temperature}, since it is used
 to prove the main results of the paper, which also consider deterministic dynamics.
\end{remark}

\section{Frequency control of thermostatic loads}\label{Sec:Novel_thermostatic}

{In this section} we present a frequency control scheme for TCLs and propose appropriate conditions for its design.
For the proposed scheme, we show that, as the number of loads tends to infinity, then  no synchronization phenomena occur and
that there exist arbitrarily long time intervals where frequency deviations are arbitrarily small.

\subsection{Frequency control scheme for thermostatic loads}
\begin{figure}[t!]
\centering
\includegraphics[scale = 0.6]{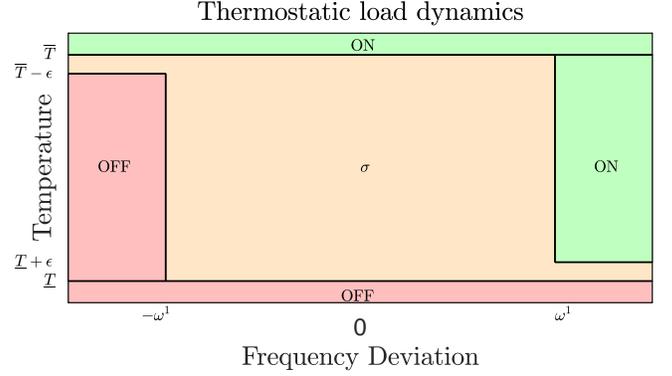}
\caption{TCL scheme described by \eqref{sys_hysteresis_3}. {In the green and red areas the switching state is ON and OFF respectively. In the orange {area}  the switching state can be either ON or OFF.}}
\label{Scheme_figure}
\end{figure}

{We introduce in this subsection the frequency control policy for the TCLs, which is a scheme that provides an ancillary service at urgencies,} i.e. when frequency deviations exceed particular {thresholds.}
 The scheme, depicted in Figure \ref{Scheme_figure}, is described below
\begin{subequations}\label{sys_hysteresis_3}
\begin{gather}
d^c_j = \overline{d}_j \sigma_{j}, \\
\hspace{-0.5mm}\sigma_{j}(t^+) \hspace{-0.5mm}= \hspace{-0.5mm}
\begin{cases}
 1, \hspace{2mm} \quad  \begin{cases} T_j \geq  \overline{T}_j,\\[1mm]
 {\omega \geq  \omega^1_j \text{ and }  {T_j \geq \underline{T}_j + \epsilon_j}},
\end{cases} \\[1mm]
0, \hspace{2mm}  \quad  \begin{cases}  T_j \leq  \underline{T}_j,\\[1mm]
 {\omega \hspace{-0.5mm}\leq \hspace{-0.5mm} -\omega^1_j \text{ and }  {T_j \hspace{-0.5mm}\leq\hspace{-0.5mm} \overline{T}_j \hspace{-0.5mm}- \hspace{-0.5mm}\epsilon_j}},
 \end{cases} \\[1mm]
\sigma_j(t),    \begin{cases}  |\omega| \leq \omega^1_j \text{ and } \underline{T}_j \leq T_j \leq \overline{T}_j,\\[1mm]
{\omega \hspace{-0.5mm}\leq \hspace{-0.5mm} -\omega^1_j \text{ and }  T_j \in {[\overline{T}_j \hspace{-0.5mm}- \hspace{-0.5mm}\epsilon_j,  \overline{T}_j]},} \\
{\omega \geq  \omega^1_j \text{ and }  T_j \in [\underline{T}_j, \underline{T}_j \hspace{-0.5mm} + \hspace{-0.5mm} \epsilon_j]},
\end{cases}
\end{cases} \label{sys_hysteresis_3_b}
\end{gather}
\end{subequations}
where $\omega^1_j > 0$ are frequency thresholds and $0 < \epsilon_j  < (\overline{T}_j - \underline{T}_j)/2, j \in N$. Note that, $\epsilon_j$  in \eqref{sys_hysteresis_3} serves to ensure than no Zeno behavior {occurs}  as a result of the coupling between the frequency and TCL dynamics. The latter is analytically shown  in Lemma \ref{dwell_time_lemma} below.

The scheme in \eqref{sys_hysteresis_3} responds to frequency deviations by switching when prescribed frequency thresholds are exceeded {thus} providing ancillary services to the power network.
 Furthermore, when the frequency deviation does not reach the corresponding frequency thresholds, then the scheme in \eqref{sys_hysteresis_3} reduces to  \eqref{sys_hysteresis}.
  Note that, according to \eqref{sys_hysteresis_3}, {the} temperature will always be within its respective thresholds and hence users comfort levels will not be affected.

For the rest of the manuscript, we let $S(\bar{\omega}) = \{{j \in N}: \omega^1_j \leq \bar{\omega} \}$ be the set of loads with respective frequency thresholds below $\bar{\omega}$.
Moreover, for any set $S \subseteq N$, we let $\omega_m(S) = \min_{j \in S} \omega^1_j, d^{s}_{S}(t) = \sum_{j \in S} d^c_j(t)$ and $d^{s,*}_{S} = \sum_{j \in S} \alpha_j \overline{d}_j$.
Furthermore, we {let $\hat{L}$ be the 1-norm of the system with input $d^s$ and output $\omega$, described by \eqref{sys1_2}, {which is given by}
\begin{gather}
\hat{L} = \int_{0}^{\infty} |C\mathrm{e}^{At}B| dt, \label{eq:1norm}
\end{gather}
where  $C = [1 \; \vect{0}^T_n]$, noting that its boundedness follows from Assumption \ref{assum_cont_system}(i).}

The following condition is imposed for the design of frequency thresholds. Within it, we let $\zeta_j = \max(\alpha_j, 1-\alpha_j)$, noting that $\zeta_j \in {[0.5},1)$ since $\alpha_j \in (0,1)$.

\begin{design condition}\label{des_condition_freq_thresholds}
{The frequency thresholds $\omega^1_j$ are chosen such that  for all $\bar{\omega} \in \mathbb{R}_+$ and some $\delta > 0$,
 $\sum_{j \in S(\bar{\omega})} \zeta_j \overline{d}_j \leq  \max({(\bar{\omega} - \delta)/\hat{L}},0)$, where $\hat{L}$ is given {by \eqref{eq:1norm}.}}
\end{design condition}

Design condition \ref{des_condition_freq_thresholds} restricts the coupling of frequency and TCL dynamics by bounding the aggregate demand that actively contributes to frequency regulation.
 The condition allows to deduce that no synchronization occurs between TCLs when the scheme \eqref{sys_hysteresis_3} is implemented.
Note also that $\delta$ in Design condition \ref{des_condition_freq_thresholds} satisfies  $\delta \in (0, \omega_m(N))$ by definition, since $\omega^1_j < \delta$ for some $j \in N$ would imply that Design condition \ref{des_condition_freq_thresholds} does not hold.
{To implement Design condition \ref{des_condition_freq_thresholds}, the values of $\omega^1$ for the TCL population should be selected such that the presented bound is satisfied at all values of $\bar{\omega}$, i.e. given  $\bar{\omega}$, the condition restricts the {aggregate demand}  
of loads that may switch due to that particular frequency deviation.}

 \subsection{Hybrid system description}

The behavior of system \eqref{sys1}, \eqref{sys2}, \eqref{sys_temperature}, \eqref{sys_hysteresis_3}, can be described by the states {$z = (\overline{x},\sigma$), where $\overline{x} =  (\omega, \hat{x}, T) \in \mathbb{R}^m$, $m = |N| + n + 1$,}
is the continuous state, and $\sigma \in P^{|N|}$ the discrete state, where $P = \{0,1\}$. Moreover, {we} let $\Lambda =  \mathbb{R}^m \times P^{|N|}$ be the space where the system states evolve.
The continuous dynamics of the system \eqref{sys1}, \eqref{sys2}, \eqref{sys_temperature}, \eqref{sys_hysteresis_3}, are described by
\begin{subequations}\label{sys4_hysteresis}
\begin{gather}
\hspace{-5mm} M \dot{\omega} = - p^L + p^M -  D \omega - \sum_{j \in N} \overline{d}_j\sigma_j,  \label{sys4a} \\
{p^M = \hat{C} \hat{x} + \hat{D} \omega,} \quad {\dot{\hat{x}} = \hat{A} \hat{x} + \hat{B} \omega,} \label{sys4b2}\\
 \dot{T}_j = -k_j (T_j - \hat{T}_j + \lambda_j {\overline{d}_j\sigma_j}), j \in N, \label{sys4c}  \\
 \dot{\sigma}_j = 0, j \in N, \label{sys4d}
\end{gather}
\end{subequations}
which is valid when $z$ belongs to the set $F$ given by
\begin{equation}\label{e:C}
F=\{ z \in \Lambda:  \sigma_j \in \mathcal{I}_j(T_j, \omega), \;\forall j\in N\},
\end{equation}
where
\[
\mathcal{I}_j(T_j, \omega)= \begin{cases}
\{1\}, \begin{cases}
 T_j > \overline{T}_j,\\[1mm]
 {\omega >  \omega^1_j \text{ and } {T_j > \underline{T}_j + \epsilon_j}},
\end{cases}
  \\[1mm]
\{0\}, \begin{cases}
 T_j < \underline{T}_j,\\[1mm]
 {\omega \hspace{-0.5mm}<\hspace{-0.5mm}  -\omega^1_j \text{ and } {T_j \hspace{-0.5mm} < \hspace{-0.5mm} \overline{T}_j - \epsilon_j}},
 \end{cases} \\[1mm]
\{0, 1\}, \begin{cases}  |\omega| \leq \omega^1_j \text{ and } \underline{T}_j \leq T_j \leq \overline{T}_j,\\[1mm]
 {\omega \hspace{-0.5mm}\leq \hspace{-0.5mm} -\omega^1_j \text{ and }  T_j \in {[\overline{T}_j \hspace{-0.5mm}- \hspace{-0.5mm}\epsilon_j,  \overline{T}_j]},} \\
{\omega \geq  \omega^1_j \text{ and }  T_j \in [\underline{T}_j, \underline{T}_j \hspace{-0.5mm} + \hspace{-0.5mm} \epsilon_j]}.
\end{cases}
\\[1mm]
\end{cases}
\]

Alternatively, when $z$ belongs to the set $G = {(\Lambda \setminus F) \cup \underline{G}}$
where
$\underline{G} = \{ z \in \Lambda:
\sigma_j \in \mathcal{I}^D_j(T_j, \omega), \;\forall j\in N\}$,
and
\[
\mathcal{I}^D_j(T_j, \omega)= \begin{cases}
\{1 \},  \begin{cases}
{\omega \geq -\omega^1_j \text{ and }} T_j =  \underline{T}_j, \\
\omega = -\omega^1_j \text{ and } {T_j \in [\underline{T}_j, \overline{T}_j \hspace{-0.5mm}-\hspace{-0.5mm} \epsilon_j}], \\
\omega \leq -\omega^1_j \text{ and } {T_j =  \overline{T}_j \hspace{-0.5mm}-\hspace{-0.5mm} \epsilon_j},
\end{cases} \\
\{0\},
\begin{cases} {\omega \leq \omega^1_j \text{ and }} T_j =  \overline{T}_j,\\[1mm]
\omega = \omega^1_j \text{ and } {T_j \in [\underline{T}_j + \epsilon_j, \overline{T}_j]}, \\
\omega \geq \omega^1_j \text{ and } {T_j = \underline{T}_j + \epsilon_j},
\end{cases}
\end{cases}
\]
then its components follow the discrete update described below
\begin{align}\label{sys4_g}
\overline{x}^+ &= {\overline{x}(t)}, \;
\sigma_{j}(t^+) = \begin{cases}
1,\hspace{-0.5mm}
\begin{cases}
T_j \geq \overline{T}_j,\\[1mm]
\hspace{-0.5mm} \omega \hspace{-0.5mm} \geq \hspace{-0.5mm}  \omega^1_j  \text{ \hspace{-0.25mm}and\hspace{-0.25mm} } T_j \hspace{-0.5mm} \in\hspace{-0.5mm} [\underline{T}_j \hspace{-0.5mm}+\hspace{-0.5mm} \epsilon_j, {\overline{T}_j}],
\end{cases} \\
0,\hspace{-0.5mm}
\begin{cases}
T_j \leq \underline{T}_j,  \\[1mm]
\hspace{-0.5mm}\omega \hspace{-0.5mm}\leq \hspace{-0.5mm} -\omega^1_j\hspace{-0.5mm} \text{ \hspace{-0.25mm}and\hspace{-0.25mm} } T_j \hspace{-0.5mm}\in\hspace{-0.5mm} [{\underline{T}_j}, \overline{T}_j \hspace{-0.75mm} -\hspace{-0.75mm} \epsilon_j],
\end{cases} \\
\end{cases}
\end{align}
{where $\overline{x}^+ = \lim_{\epsilon \rightarrow 0} \overline{x}(t + \epsilon)$.}

 We can now provide the following compact representation for the hybrid system \eqref{sys1}, \eqref{sys2}, \eqref{sys_temperature}, \eqref{sys_hysteresis_3},
\begin{equation}\label{sys4}
\dot{z} = f(z), z \in F, \quad
z^+ = g(z), z \in G,
\end{equation}
where $f(z): F \rightarrow \Lambda$ and $g(z): G \rightarrow F$  are described by \eqref{sys4_hysteresis} and \eqref{sys4_g} respectively.  Note that $z^+ = g(z)$ represents a discrete dynamical system where $z^+$ indicates that the next value of the state $z$ is given as a function of its current value through $g(z)$. Moreover, notice that $F \cup G = \Lambda$.

\subsection{Analysis of solutions}

In this section we consider the solutions of \eqref{sys4} and show their existence and that no Zeno behavior occurs.
Below we provide a definition of a hybrid time domain, hybrid solution and complete and maximal solutions for systems described by \eqref{sys4} from \cite[Ch. 2]{goebel2012hybrid}. \ka{Note that the definition of a hybrid system is provided in \cite[Dfn. 2.2]{goebel2012hybrid}.}

\begin{definition}(\cite{goebel2012hybrid})\label{dfn_hybrid_domain_solution}
A subset of $\mathbb{R}_{\geq 0} \times \mathbb{N}_0$ is a hybrid time domain if it is a union of a finite or infinite sequence of intervals  $[t_\ell, t_{\ell+1}] \times \{\ell\}$, with the last interval (if existent) possibly of the form $[t_\ell, t_{\ell+1}] \times \{\ell\}$, $[t_\ell, t_{\ell+1}) \times \{\ell\}$, or $[t_\ell, \infty) \times \{\ell\}$.
Consider a function $z(t,\ell): K \rightarrow \mathbb{R}^m$ defined on a hybrid time domain $K$ such that for every fixed $\ell \in \mathbb{N}$, $t\rightarrow z(t,\ell)$ is locally absolutely continuous on the interval $T_\ell = \{t: (t,\ell) \in K\}$.
 The function $z(t,\ell)$ is a solution to the hybrid system $\mathcal{H} = (F,f,G,g)$ if {$z(0,0) \in {{F} \cup G}$}, and for all $\ell \in \mathbb{N}$ such that $T_\ell$ has non-empty interior (denoted by ${\rm int}T_l$)
\begin{align*}
& z(t,\ell) \in F, \text{ for all t} \in {{\rm int}T_l}, \\
& \dot{z}(t,\ell) \in f(z(t,\ell)), \text{ for almost all } t \in T_\ell,\\
& \hspace{-12mm} \text{and for all } (t,\ell) \in K \text{ such that } (t,\ell+1) \in K,\\
& z(t,\ell) \in G, \; z(t,\ell+1) \in g(z(t,\ell)).
\end{align*}
A solution $z(t,\ell)$ is complete if $K$ is unbounded. A solution $z$ is maximal if there does not exist another solution $\tilde z$ with time domain $\tilde K$ such that $K$ is a proper subset of $\tilde K$ and $z(t,j) = \tilde{z}(t,j)$ for all $(t,j) \in K$.
\end{definition}

The following lemma, proven in the appendix, {shows the existence  of complete solutions to \eqref{sys4}. Furthermore, it demonstrates the boundedness of solutions to \eqref{sys4} and provides a lower bound on the time between  consecutive switches, which suffices to show that no Zeno behavior occurs.
Finally, it states that all maximal solutions to \eqref{sys4} are complete.}
We remind that $t_{j,i}$ is the time of the $i$th switch of load $j$.

{\begin{lemma}\label{dwell_time_lemma}
For any initial condition $z(0,0) \in {\Lambda}$ there exists a complete solution to \eqref{sys4}.
Furthermore, all maximal solutions to \eqref{sys4} are complete.
Moreover, if  Assumption \ref{assum_cont_system} holds then the following hold:
\begin{enumerate}[(i)]
\item  {For each initial condition $z(0,0) \in {\Lambda}$, solutions to \eqref{sys4}  are bounded.}
\item For any
 solution to \eqref{sys4}, there exists ${\tau_d} > 0$  such that $\min_{i \geq 1} (t_{j, i+1} - t_{j,i}) \geq {\tau_d}$ for any $j \in N$.
\end{enumerate}
\end{lemma}}
 {The boundedness of solutions to \eqref{sys4}, demonstrated in the above lemma, follows also intuitively by noting that \eqref{sys4} consists of the asymptotically stable linear system \eqref{sys1}, \eqref{sys2}, {with input $d^s$ and output $\omega$} in feedback with the hybrid system \eqref{sys_temperature}, \eqref{sys_hysteresis_3} and that the magnitude of {$d^s$, which can be regarded as the output of \eqref{sys_temperature}, \eqref{sys_hysteresis_3},} is bounded.
 Furthermore, the boundedness of $T_j, j \in N$ follows directly from the structure of \eqref{sys_temperature}, \eqref{sys_hysteresis_3}.}

\subsection{Performance analysis}

In this section we state {one of the main results} of this paper, {associated with} the  performance of solutions to \eqref{sys4}.
The following theorem, proven in  the appendix, demonstrates that as the number of loads tends to infinity, then for all initial conditions
there exist arbitrarily long time intervals where frequency deviations are arbitrarily small.

\begin{theorem}\label{thm_conv}
Consider the system described by \eqref{sys4} and let Assumptions {\ref{assum_cont_system}--\ref{assumption_period_ratio}} and Design condition \ref{des_condition_freq_thresholds} hold. Furthermore, assume that the thermostatic loads  described by \eqref{sys_temperature}, \eqref{sys_hysteresis_3} satisfy $\overline{d}_j = \frac{\Gamma}{|N|}$. Then, as $|N| \rightarrow \infty$, for  any $z(0,0)~\in~\mathbb{R}^m~\times~P^{|N|}$, any maximal solution of \eqref{sys4} and any $\epsilon > 0, \hat{\tau} \in \mathbb{R}_+$, there exists $\tau \in \mathbb{R}_+$   such that
 $|\omega(t,j)| \leq \epsilon$ for $t \in [\tau, \tau + \hat{\tau}]$.
\end{theorem}

{The importance of Theorem \ref{thm_conv} is that it
  shows, {for all initial conditions}, that frequency trajectories become  arbitrarily small for an arbitrarily long amount of time.
Also, as shown in {Lemma \ref{dwell_time_lemma}} the scheme in \eqref{sys_hysteresis_3} avoids Zeno behavior. {Furthermore,} being deterministic,
it allows the instant response to  frequency deviations, thus providing {improved} ancillary services to the power system. The latter, is {also} numerically demonstrated in the following section.}

\begin{remark}\label{remark_tau}
Theorem \ref{thm_conv} does not provide an analytical expression for $\tau$.
However, it is intuitive to note that its value in a real setting depends on: (i) the values of $\epsilon$ and $\hat{\tau}$, which are associated with its definition, (ii) the initial conditions and the speed of generation dynamics, which determine how long it takes for generation to match a potential disturbance,
and (iii) the distribution of load {periods.}
\end{remark}

\section{Simulation on the NPCC 140-bus system} \label{Simulation}

In this section we verify our analytic results with a numerical simulation on the Northeast Power Coordinating
Council (NPCC) 140-bus interconnection system, using the Power System Toolbox~\cite{cheung2009power}. This model is more detailed and realistic than our analytical one, including line resistances, a DC12 exciter model{, a transient reactance generator model, and  turbine governor dynamics.

The test system consists of $93$ load buses serving different types of loads including constant active and reactive loads and $47$ generation buses. The overall system has a total real power of $28.55$~GW. For our simulation, we added five loads on buses $2, 8, 9, 16$ and $17$, each having a step increase of magnitude $2$ p.u. (base $100$MVA) at $t=1$ second.

Controllable loads were considered within the simulations at load buses $1 - 20$, with loads controlled every $10$ms. In particular, we considered $500$ refrigerators of equal magnitude at each of  the $20$ selected load buses with 
{aggregate}
power of\footnote{{A more realistic simulation would involve $10^6$ refrigerators for the same aggregate demand but would be computationally expensive. The simulated number suffices to demonstrate the analysis {in the paper} noting that a larger number {of TCLs} would result {in} an even  smoother~response.}}
 $2.5$~GW.
 For comparison, we considered  the system response  when the following {four} schemes for TCLs {were}~implemented.
 \begin{enumerate}[(i)]
 \item   Conventional {TCLs} 
 {that do not contribute to frequency control}, i.e. loads with dynamics as in \eqref{sys_hysteresis}, \eqref{sys_temperature}.
 \item  {Frequency dependent {TCLs}  
 with a deterministic control policy}, i.e. loads with dynamics described by \eqref{sys_temperature}, \eqref{sys_hysteresis_3}.
\item   Frequency dependent {TCLs}  {with a randomized control policy},
 {as in} \cite{angeli2012stochastic}, \cite{aunedi2013economic}.
 \item The scheme (iii) with larger feedback gains, aiming for a faster response.
 \end{enumerate}
 The above cases  will be referred {to as} {case (i), (ii), (iii) and {(iv)}} respectively.
 The {values of the {control} parameters} were randomly selected from uniform distributions with bounds provided {in Table \ref{table_values}}.
 Furthermore, initial conditions were randomly selected in a similar manner.
 To ensure that incorporating the loads would not disturb the balance of the network,
   for each thermostatic load {incorporated at a bus}  some constant demand equal to its average value was removed from the same bus.
  Moreover, frequency thresholds in case {(ii)} were selected in accordance with Design condition \ref{des_condition_freq_thresholds}.
{In particular, following the approach described in \cite{shi2018analytical}, an equivalent single bus model of the power network, where generation was described with high order dynamics, was derived.
The latter enabled to obtain $\hat{L}$ (i.e. the $1$-norm of the system \eqref{sys1_2} with {input the aggregate demand and output $\omega$) and implement Design condition \ref{des_condition_freq_thresholds}.
To ensure that Design condition \ref{des_condition_freq_thresholds} was satisfied, we {verified that the
 selected} values of $\omega^1$ satisfied $\sum_{j \in S(\bar{\omega})} \zeta_j \overline{d}_j \leq  \max((\bar{\omega} - \delta)/\hat{L},0)$, letting $\delta = 0.001Hz$, for all $\bar{\omega} \in {\mathbb{R}_+}$.}
For additional safety,  frequency thresholds {were} designed with a $20\%$ margin from the obtained upper~bound.}
 For case (iii), the implemented algorithm involved {randomized transitions between the on/off states with controlled rates\footnote{{Additional temperature constraints were not considered for simplicity, as simulations indicate that these restrict the frequency {control} performance when the temperature thresholds are fixed.}} as in \cite{angeli2012stochastic}.}
  The algorithm was implemented with {$K_{\pi} = 5$ and  $v_{des} = 1$ for each TCL} (in analogy to \cite{angeli2012stochastic}), where  {$K_{\pi}$ and $v_{des}$} {are parameters associated with
  the  {feeback gain} and the desired temperature variability} respectively.
  For case (iv), we implemented case (iii) with  $K_{\pi} = 50$.

\begin{table}[t!]
\centering
 \begin{tabular}
 {|p{0.07\textwidth} | p{0.1\textwidth} | p{0.1\textwidth} |}
 \hline
  Variable & Lower Bound & Upper Bound \\ [0.1ex]
  \hline
   \hline
  $\hat{T}$ & \SI{15}{\celsius} & \SI{25}{\celsius} \\  [0.1ex]
 \hline
  $\overline{T}$ & \SI{5}{\celsius} & \SI{7}{\celsius} \\
 \hline
   $\underline{T}$ & \SI{2}{\celsius}  & \SI{4}{\celsius}  \\
 \hline
    $k$ & 2$\times 10^{-4}$ & $10^{-3}$ \\
 \hline
     $\lambda_j$ & $25 (\overline{d}_j)^{-1}$ & $35 (\overline{d}_j)^{-1}$ \\
 \hline
     $\omega^1$ & 0.01 $Hz$ & 0.26 $Hz$ \\
 \hline
     $\epsilon_j$ & \SI{0.001}{\celsius} & \SI{0.01}{\celsius} \\
 \hline
\end{tabular}
\vspace{-0mm}
   \caption{Ranges of coefficients describing TCL dynamics.}
   \label{table_values}
\end{table}

The frequency at bus 27 for the {four} 
tested cases is  shown in Figure \ref{Frequency}.
{We} observe that {the} frequency converges to a very small set containing its nominal value.
Furthermore, Figure \ref{overshoot_freq} suggests that  the scheme in \eqref{sys_hysteresis_3} results in {a reduced} frequency overshoot {relative to the other cases considered, by {illustrating} the largest deviation in frequency}
at buses $1-40$,
{which are {the} buses} where {the frequency} overshoot was seen to be the largest.
{In addition, it demonstrates that increasing the feedback gains  {in the schemes with randomization} {in cases (iii), (iv)} results in {a} reduced frequency overshoot. However, 
{larger transition rates can lead to  more frequent switching of the TCLs, which is generally undesirable. Furthermore, increasing only one of the transition rates (as is the effect of increasing {$K_{\pi}$}) will maintain a slower recovery of the TCLs. The speed of response of schemes {with randomization} can potentially be improved by combining them with deterministic schemes, which is an interesting direction for further theoretical analysis.
Figures \ref{Frequency}, \ref{overshoot_freq} also demonstrate that 
no Zeno behavior or load synchronization are experienced with the proposed deterministic scheme}.}

The percentage of TCLs that are ON  for {the TCL schemes} described in cases (i) {and} (ii) is depicted in Figure \ref{Fig_thm_1}. {It should be noted that the almost flat response in case (i) validates Theorem \ref{thm_variance}.}

\begin{figure}[t!]
\centering
\includegraphics[scale = 0.63]{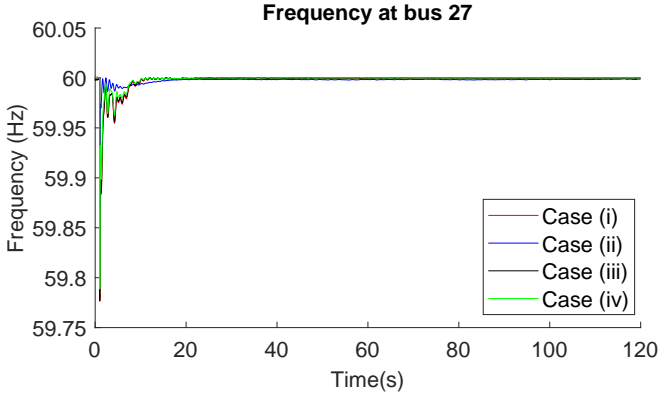}
\vspace{-1mm}
\caption{Frequency at bus 27 with TCL dynamics {in the} following cases: i) Conventional {TCLs},
ii)~Deterministic frequency dependent {TCLs,} 
iii)~Frequency dependent {TCLs} {with randomized control policy},
iv) As in (iii) but feedback gains are ten times larger.}
\label{Frequency}
\end{figure}
\begin{figure}[t!]
\centering
\includegraphics[scale = 0.63]{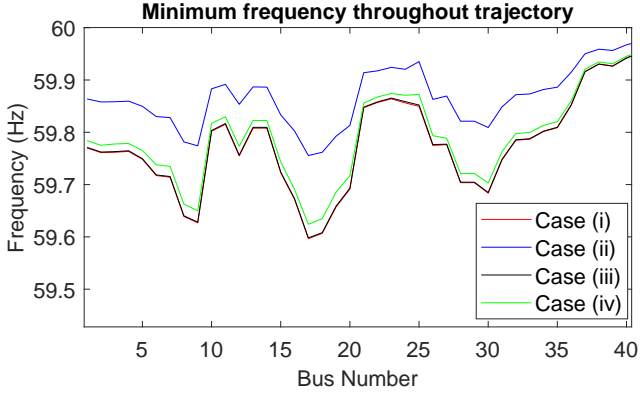}
\vspace{-1mm}
\caption{Largest frequency overshoot for buses $1-40$ for the {four cases described in the caption of {Figure}~\ref{Frequency}.}}
\label{overshoot_freq}
\end{figure}
\begin{figure}[t!]
\centering
\includegraphics[scale = 0.63]{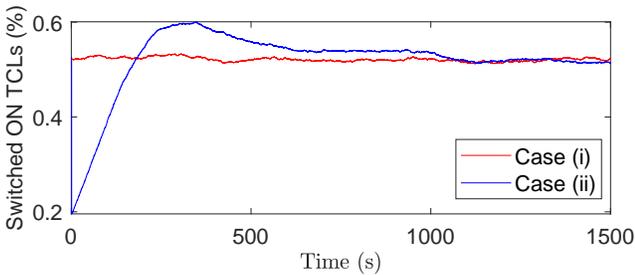}
\vspace{-1mm}
\caption{{Percentage of TCLs switched ON for the following cases: i) Conventional {TCLs,} 
ii)~Deterministic frequency dependent {TCLs.}}} 
\label{Fig_thm_1}
\end{figure}

\section{Conclusion}\label{Conclusion}

We have {studied} the problem of controlling thermostatic loads to provide ancillary services to the power network at urgencies.
We first considered conventional TCLs and showed that their aggregation has zero variance when their number tends to infinity and a 
{mild} condition on their period ratios holds.
Then, we proposed a deterministic control scheme for TCLs which induces switching when frequency deviations exceed particular frequency thresholds.
For the considered scheme, we explain how frequency thresholds {could} be designed such that the coupling between load and frequency dynamics does not cause load synchronization.
{In particular,} when the number of loads tends to infinity, we showed that frequency deviations are arbitrarily small for arbitrarily large periods of {time.}
Our analytic results have been numerically verified with simulations on the NPCC 140-bus system, which demonstrate improved frequency response when frequency dependent TCLs are incorporated  compared to when conventional implementations {are} considered.

{Future extensions of this work could consider more involved dynamics, including a network model of the power grid.
In addition, future studies could consider more advanced control designs for TCLs, taking into account elements such as their economic performance and the rate at which they desynchronize, which may yield improved response.}

\appendix


\emph{Proof of Theorem \ref{thm_variance}:} By definition, \ilcc{the} variance is given by
\begin{equation}\label{var_eqt_defn}
\mathbb{V}(d^s) = \mathbb{E}((d^s)^2) - [\mathbb{E}(d^s)]^2.
\end{equation}
Since $d^s = \sum_{j \in N} d^c_j$, it then holds that
\begin{align}\label{var_first_term}
\mathbb{E}((d^s)^2) &= \sum_{j \in N} \mathbb{E}((d^c_j)^2) + 2\sum_{(i,j) \in E} \mathbb{E}(d^c_id^c_j) \nonumber \\
&= \sum_{j \in N} \alpha_j (\frac{\Gamma}{|N|})^2 + 2\sum_{(i,j) \in E} \alpha_i \alpha_j  (\frac{\Gamma}{|N|})^2,
\end{align}
where the first argument follows trivially and the second from Proposition \ref{prop_cross_term} below.
Furthermore, the second term in \eqref{var_eqt_defn} satisfies
\begin{align}\label{var_second_term}
&[\mathbb{E}(d^s)]^2 = [\mathbb{E}(\sum_{j \in N} d^c_j)]^2  = \sum_{j \in N} \mathbb{E}(d^c_j)^2 + 2\hspace{-0.5mm}\sum_{(i,j) \in E}\hspace{-0.5mm} \mathbb{E}(d^c_i) \mathbb{E}(d^c_j) \hspace{-0.5mm} \nonumber \\
&= \hspace{-0.5mm} \sum_{j \in N}\hspace{-0.5mm} \alpha_j^2 (\frac{\Gamma}{|N|})^2 + 2\hspace{-0.5mm}\sum_{(i,j) \in E}\hspace{-0.5mm} \alpha_i \alpha_j  (\frac{\Gamma}{|N|})^2.
\end{align}

Combining \eqref{var_eqt_defn}, \eqref{var_first_term} and \eqref{var_second_term} results to
$
\mathbb{V}(d^s) = \sum_{j \in N} \alpha_j (1 - \alpha_j) (\frac{\Gamma}{|N|})^2 \antk{<} \frac{\Gamma^2}{|N|}
$
noting for the last argument that $0 < \alpha_j < 1, j \in N$. Hence, it holds that $\lim_{|N| \rightarrow \infty} \mathbb{V}(d^s) = 0$.
 \hfill $\blacksquare$

\begin{proposition}\label{prop_cross_term}
 Consider TCLs  described by \eqref{sys_hysteresis}, \eqref{sys_temperature} and let Assumption \ref{assumption_period_ratio} hold. Then,  $\mathbb{E}(d^c_id^c_j)~=~\alpha_i \alpha_j \Gamma_i \Gamma_j$ for all $(i,j) \in E$.
\end{proposition}

\emph{Proof of Proposition \ref{prop_cross_term}:}
From Assumption \ref{assumption_period_ratio} it follows that  $\rho_{ij} \in \mathbb{R}_+ / \mathbb{Q}_+$ for all $(i,j) \in E$ and hence the signal $d^c_i(t)d^c_j(t)$ is aperiodic. Its average is defined as
\begin{equation}\label{mean_dfn_aperiodic}
\mathbb{E}(d^c_id^c_j) = \lim_{\tau \rightarrow \infty} \frac{1}{\tau} \int_{0}^{\tau} d^c_i (t) d^c_j(t) dt.
\end{equation}

Without \ilcc{loss of generality} 
let $\pi_i > \pi_j$ and $t_{i,1}, t_{j,1}$ be the first time when loads $i$ and $j$ switch ON respectively.
Then,  let
 $\hat{t}_k = [t_{j,1} - t_{i,1} - (k-1)\pi_i]^+_{\pi_j}$  noting that it represents the time difference between the $k$-th time load $i$ switches ON and the first time load $j$ switches ON afterwards.
Then, letting\footnote{\ank{The fact that $[\hat{t}_{k+1} - \hat{t}_k]^+_{\pi_j} \ank{= [-\pi_i]^+_{\pi_j}}$ follows from the modular addition property, which is a standard property in modular arithmetics, see e.g. \cite[Ch. 4.1]{rosen2012discrete}.}}
 $c_k = [\hat{t}_{k+1} - \hat{t}_k]^+_{\pi_j} \ank{= [-\pi_i]^+_{\pi_j}}$, it follows that $c_k = c, \forall k \geq 1$.
Furthermore, \ank{since $c = \mu \pi_j - \pi_i$, for some $\mu \in \mathbb{N}$ and} $\rho_{ij} \in \mathbb{R}_+ / \mathbb{Q}_+$, it follows that $\frac{c}{\pi_j} \in \mathbb{R}_+ / \mathbb{Q}_+$.
 Hence, \eqref{mean_dfn_aperiodic} satisfies
\begin{equation}\label{mean_dfn_aperiodic_ext}
\mathbb{E}(d^c_id^c_j) = \frac{1}{\pi_i} \lim_{\hat{N} \rightarrow \infty} \frac{1}{\hat{N}} \sum_{k = 1}^{\hat{N}} \int_{0}^{\pi_i} d^c_i (t) d^c_j(t - \hat{t}_k) dt.
\end{equation}

From the definition of $\hat{t}_k$ it follows that its values lie within $[0, \pi_j]$. Furthermore, from $c_k = c, k \geq 1$ and $\frac{c}{\pi_j} \in \mathbb{R}_+ / \mathbb{Q}_+$, it follows that the \ank{sequence of $[\hat{t}_k]^+_{\pi_j}$ becomes uniformly distributed} as $k \rightarrow \infty$.
The latter follows by noting that the sequence \ank{$\{[\hat{t}_2 - \hat{t}_1]^+_{\pi_j}, [\hat{t}_3 - \hat{t}_1]^+_{\pi_j}, [\hat{t}_4 - \hat{t}_1]^+_{\pi_j}, \dots, [\hat{t}_{N+1} - \hat{t}_1]^+_{\pi_j}, \dots \}$, which is equal to $\{[c]^+_{\pi_j}, [2c]^+_{\pi_j}, [3c]^+_{\pi_j}, \dots, [Nc]^+_{\pi_j}, \dots \}$,}
is equivalent to the sequence $\pi_j \times \{[\frac{c}{\pi_j}]^+_{1}, [\frac{2c}{\pi_j}]^+_{1}, [\frac{3c}{\pi_j}]^+_{1}, \dots, [\frac{Nc}{\pi_j}]^+_{1}, \dots \}$ and the uniformity of this sequence is
 a special case of the Weyl Criterion (e.g. \cite[Theorem 2.1]{kuipers2012uniform}) since $\frac{c}{\pi_j}$ is irrational.
From the last argument it follows that \eqref{mean_dfn_aperiodic_ext} can be equivalently written as
\begin{equation*}
\mathbb{E}(d^c_id^c_j) = \frac{1}{\pi_i \pi_j} \int_{0}^{\pi_j} \int_{0}^{\pi_i} d^c_i (t) d^c_j(t-\bar{t}_1) dt d\bar{t}_1.
\end{equation*}
Then, considering that for  $t \in [0, \pi_i]$ it holds that
\[
d^c_i(t)
 =
 \begin{cases} \Gamma_i, t \in [t_{i,1}, t_{i,1} + \alpha_i \pi_i],  \\
 0, \text{otherwise},
 \end{cases}
\]
 and defining $s = \max_{\gamma \in \mathbb{N}_0} \{ \gamma : \alpha_i \pi_i > \gamma \pi_j \}$,  it follows that
\begin{multline}\label{second_int_part}
\mathbb{E}(d^c_id^c_j) = \frac{\Gamma_i}{\pi_i \pi_j} \int_{0}^{\pi_j} [s\pi_j \alpha_j \Gamma_j
+      \int_{t_{i,1} + s\pi_j}^{t_{i,1} + \alpha_i \pi_i} \hspace{-4mm} d^c_j(t-\bar{t}_1) dt] d\bar{t}_1.
\end{multline}
The second integral in \eqref{second_int_part} can be evaluated as
\begin{equation*}
\int_{0}^{\pi_j} \int_{t_{i,1} + s\pi_j}^{t_{i,1} + \alpha_i \pi_i}  d^c_j(t-\bar{t}_1) dt d\bar{t}_1 = \ank{\pi_j} (\alpha_i \pi_i - s\pi_j) \alpha_j \Gamma_j,
\end{equation*}
which from \eqref{second_int_part} results to $\mathbb{E}(d^c_id^c_j) = \alpha_i \alpha_j \Gamma_i \Gamma_j$.
 \hfill $\blacksquare$

\an{
\emph{Proof of Lemma \ref{dwell_time_lemma}:}
The existence of a complete solution to \eqref{sys4} follows trivially from the fact that the dynamics in \eqref{sys4_hysteresis} are globally Lipschitz and that $f$ and $g$ map into $\Lambda$ which is the domain of \eqref{sys4}.
Furthermore, the fact that all maximal solutions to \eqref{sys4} are complete follows from the global Lipschitz property of $f$
 and the fact that $f$ and $g$ map into $\Lambda$ which is the domain of \eqref{sys4}  \cite[Proposition 6.10]{goebel2012hybrid}.
\ill{The rest two parts of the Lemma are proved below:}
\begin{enumerate}[(i)]
\item \ank{
\ilcc{The} boundedness of $(\omega, p^M)$ and $T_j, j \in N$ follows since \eqref{sys4a}--\eqref{sys4b2} and \eqref{sys4c} \illl{can be seen as} asymptotically stable linear systems 
(\ilcc{Assumption} \ref{assum_cont_system}(i)) with bounded inputs $\Sigma_{j \in N} d^c_j\sigma_j$ and \illl{$d^c_j\sigma_j$} respectively.}
\item \ilcc{From (i) note} that \kas{for each $z(0,0) \in \Lambda$, the} solution to system~\eqref{sys4}, with states $z = (\omega, p^M, T, \sigma)$, is bounded.
Then,
note that
 the values of $\dot{\omega}$ and $\dot{T}_j$ are bounded from above by constants, $d\omega^{\max}$ and $dT^{max}_j$, as a result of the boundedness of solutions and the fact that the vector field in \eqref{sys4_hysteresis} is globally Lipschitz.
Hence, it follows that $t_{j,\ell+1}-t_{j,\ell} \geq \min(2\omega^1_j / d\omega^{\max}_j,  \epsilon_j/dT^{max}_j) = \tau_j$. Finally, let $\ant{\tau_d} = \min_{j \in N} \tau_j$ to conclude the proof.
 \hfill $\blacksquare$
\end{enumerate}
 }

The following results will be used \ill{within the} 
 proof of Theorem \ref{thm_conv}.

\begin{corollary}\label{corollary_variance}
Let Assumption \ref{assumption_period_ratio} hold and consider TCLs described by \eqref{sys_temperature}, \eqref{sys_hysteresis_3} with $\overline{d}_j = \frac{\Gamma}{|N|}$ and any set $S \subseteq N$. Then, if there exists $\tau \geq 0$ such that $|\omega(t)| < \omega_m(S)$ for all $t \geq \tau$ then $\mathbb{V}(d^s_{S}) \rightarrow 0$ as $|N| \rightarrow \infty$.
\end{corollary}

\emph{Proof of Corollary \ref{corollary_variance}:}
When for some finite $\tau$ it holds that $|\omega(t)| < \omega^1_j, t \geq \tau, j \in S$, the scheme in \eqref{sys_hysteresis_3} reduces to \eqref{sys_hysteresis} for $j \in S$.
\ank{If $|S| \rightarrow \infty$ as $|N| \rightarrow \infty$},  the proof follows directly from Theorem \ref{thm_variance} and the boundedness of $d^s_S$.
 Alternatively, \ank{if $|S| < \infty$ as $|N| \rightarrow \infty$,} then the proof follows trivially by noting that $\lim_{|N| \rightarrow \infty} \sum_{j \in S} d^c_j \leq \lim_{|N| \rightarrow \infty} |S| \frac{\Gamma}{|N|} = 0$.
 \hfill $\blacksquare$


\begin{lemma}\label{prop_dc_sum_convergence}
Consider TCLs described by  \eqref{sys_temperature}, \eqref{sys_hysteresis_3}, with $\overline{d}_j = \frac{\Gamma}{|N|}$, any set $S \subseteq N$  and let Assumption \ref{assumption_period_ratio} hold.
\ka{Then, when $|N| \rightarrow \infty$, for
any initial condition $(T(0), \sigma(0)) \in \mathbb{R}^{|N|} \times P^{|N|}$ and any $\epsilon > 0, \hat{\tau}_1 \in \mathbb{R}_+$,
there exist  $\tau, \tau_1, \hat{\tau}_0 \in \mathbb{R}_+, \tau \leq \tau_1, \tau_1 + \hat{\tau}_1 \leq \hat{\tau}_0$  such that
if $|\omega(t)| < \omega_m(S)$ for $t \in  [\tau, \hat{\tau}_0]$,
then
 $\int_{\tau_1}^{\tau_1 + \hat{\tau}_1} (d^s_S(t) - d^{s,*}_S)^2 dt \leq \epsilon$.}

\end{lemma}

\emph{Proof of Lemma \ref{prop_dc_sum_convergence}:}
From Theorem \ref{thm_variance} it follows that when Assumption \ref{assumption_period_ratio} holds for TCLs described by \eqref{sys_hysteresis}, \eqref{sys_temperature}, with $\overline{d}^c_j = \frac{\Gamma}{|N|}$ then $\lim_{|N| \rightarrow \infty} \mathbb{V}(d^s) = 0$.
Corollary \ref{corollary_variance}  extends this result to any set $S \subseteq N$, i.e. $\lim_{|N| \rightarrow \infty} \mathbb{V}(d^s_S) = 0$.
 The latter suggests that $\lim_{|N| \rightarrow \infty} \lim_{\tau \rightarrow \infty} \frac{1}{\tau} \int_{0}^{\tau} (d^s_S(t) - d^{s,*}_S)^2 dt = 0$ which follows since
\begin{align}
 & \lim_{\tau \rightarrow \infty}\hspace{-0.5mm} \frac{1}{\tau} \hspace{-0.75mm}\int_{0}^{\tau}\hspace{-0.75mm} (d^s_S(t)\hspace{-0.5mm} -\hspace{-0.5mm} d^{s,*}_S)^2 dt
   \hspace{-0.5mm}=\hspace{-0.5mm}  \lim_{\tau \rightarrow \infty}\hspace{-0.5mm} \frac{1}{\tau} \hspace{-0.75mm}\int_{0}^{\tau}\hspace{-0.75mm} (d^s_S(t))^2\hspace{-0.5mm} -\hspace{-0.5mm}  (d_S^{s,*})^2 dt \nonumber \\[-0.0mm]
  & = \lim_{\tau \rightarrow \infty} \frac{1}{\tau} \int_{0}^{\tau} (d_S^s(t))^2  dt -  (d_S^{s,*})^2 \equiv  \mathbb{V}(d_S^s),
\end{align}
where the first step follows by expanding the squared term and using the definition of $d^{s,*}_S$.

Now consider the condition on the lemma statement, and temporarily assume that
$|\omega(t)| < \omega_m(S)$ for all $t \geq \tau$. Therefore, \eqref{sys_temperature}, \eqref{sys_hysteresis_3} reduces to \eqref{sys_hysteresis}, \eqref{sys_temperature} for $t \geq \tau$.

The next part of the proof follows by contradiction. In particular, assume there exist  $\epsilon > 0$ and $\tau_1, \hat{\tau}_1 \in \mathbb{R}_+$ such that $ \int_{\overline{\tau}}^{\overline{\tau} + \hat{\tau}_1} (d^s_S(t) - d_S^{s,*})^2 dt \geq \epsilon$ for all $\overline{\tau} \in \{\tau_1 + k \hat{\tau}_1: k \in \mathbb{N}_0\}$.
Then, $\mathbb{V}(d^s_S) \geq \frac{\epsilon}{\hat{\tau}_1}$, which contradicts the result of Theorem \ref{thm_variance}.
Hence, if $|\omega(t)| < \omega_m(S)$ for all $t \geq \tau$, then for any $\epsilon > 0, \hat{\tau}_1 \in \mathbb{R}_+$,  there exists finite $\tau_1$ such that $\int_{\tau_1}^{\tau_1 + \hat{\tau}_1} (d^s_S(t) - d_S^{s,*})^2 dt \leq \epsilon$.

To conclude the proof, note that \ank{the trajectory of $d^s(t)$ depends only on the initial conditions and the trajectory of $\omega(t)$. Hence, the trajectory of $\omega(t)$ for $t \geq  \tau_1 + \hat{\tau}_1$ does not affect the fact the result that $\int_{\tau_1}^{\tau_1 + \hat{\tau}_1} (d^s_S(t) - d_S^{s,*})^2 dt \leq \epsilon$. Therefore, the condition on $\omega(t)$ reduces to  $|\omega(t)| < \omega_m(S)$ for $t \in  [\tau, \hat{\tau}_0]$, for any $\hat{\tau}_0 \geq  \tau_1 + \hat{\tau}_1$.}
%
 \hfill $\blacksquare$


Before continuing with the rest of the results, it will be convenient to note that system \eqref{sys4} consists of the linear system \eqref{sys1_2} in feedback with the hybrid system \eqref{sys_temperature}, \eqref{sys_hysteresis_3}.
\kas{\icc{Let $\hat{x}^*$ be the equilibrium value of $\hat{x}$ in \eqref{sys1_2} when $d^{s}=d^{s,*} = \sum_{j \in N} \alpha_j \overline{d}_j$. System \eqref{sys1_2}, can be equivalently written  in terms of deviations from these equilibrium values} 
as follows
\begin{equation}\label{sys_linear_description}
y  = C x, \quad \dot{x}
 =  A x + B u,
\end{equation}
where $x = \begin{bmatrix}
{\omega} \\
\hat{x}- \hat{x}^*
\end{bmatrix}$,
\icc{$u = [d^s - d^{s,*}]$,} \kas{$y = \omega$,} 
 $C = [1 \; \vect{0}^T_n]$ and $A$ and $B$ as given in the description immediately after \eqref{sys1_2}.
Furthermore, note that $A$ is Hurwitz from Assumption \ref{assum_cont_system}(i).}


\begin{lemma}\label{lemma_g2}
Consider the system \eqref{sys_linear_description}. \ic{Let $|u(t)|$ \icc{be uniformly bounded for $t\geq0$ and satisfy the property that for any $\epsilon > 0, \hat{\tau}_1 \in \mathbb{R}_+$, there exists $\tau_1 \in \mathbb{R}_+$ such that $\int_{\tau_1}^{\tau_1 + \hat{\tau}_1} u^2 dt \leq \epsilon$.} 
Then, for  any $x(0)~\in~\mathbb{R}^{n+1}$ and any $\hat{\epsilon}> 0$\kas{, $\overline{\tau}_1 \in \mathbb{R}_+$,} there exists $\tilde{\tau}_1 \in \mathbb{R}_+$ 
such that $\kas{|y(t)|} \leq  \hat{\epsilon}$ for all $t \in [\tilde{\tau}_1, \tilde{\tau}_1 + \overline{\tau}_1]$.}
\end{lemma}

\emph{Proof of Lemma \ref{lemma_g2}:}
\ank{By assumption, given $\overline{\tau}_1$, and any $\epsilon > 0$, there exists $\hat{\tau}_1 > \overline{\tau}_1$ for which there exists $\tau_1$ such that $\int_{\tau_1}^{\tau_1 + \hat{\tau}_1} u^2 dt \leq \epsilon$}.
\kas{Furthermore, $x(\tau_1)$ is \kas{uniformly bounded for $\tau_1 \geq 0$}
since \eqref{sys_linear_description} is asymptotically stable} \icc{
and $|u(t)|$ is also uniformly bounded.} 
The trajectory of $\kas{y(t)}$ for $t \in  [\tau_1, \tau_1 + \hat{\tau}_1]$ satisfies
\begin{equation}\label{omega_trajectory_tau}
|\kas{y(t)}| \leq |\kas{C \mathrm{e}^{A(t-\tau_1)} x(\tau_1)}| + |\int_{\tau_1}^t C \mathrm{e}^{A \hat{t}} B u(t - \hat{t}) d\hat{t}|.
\end{equation}
Moreover, the integral part in \eqref{omega_trajectory_tau} satisfies
\begin{align*}
&|\hspace{-0.5mm}\int_{\tau_1}^t \hspace{-1mm} C \mathrm{e}^{A \hat{t}} B u(t \hspace{-0.5mm}-\hspace{-0.5mm} \hat{t}) d\hat{t}|
\hspace{-0.5mm}\leq \hspace{-0.5mm}(\int_{\tau_1}^t \hspace{-1mm}(C \mathrm{e}^{A \hat{t}} B)^2 d\hat{t})^{\frac{1}{2}} (\hspace{-1mm}\int_{\tau_1}^t \hspace{-1mm}(u(t - \hat{t})^2 d\hat{t})^{\frac{1}{2}}  \nonumber \\
 &\leq \sqrt{\epsilon}(\int_{\tau_1}^t (C \mathrm{e}^{A \hat{t}} B)^2 d\hat{t})^{\frac{1}{2}}
\leq \sqrt{\epsilon}(\int_{0}^\infty (C \mathrm{e}^{A \hat{t}} B)^2 d\hat{t})^{\frac{1}{2}} \nonumber \\
 &=  \sqrt{\epsilon} \|\hat{G}\|_2 \leq \delta,
\end{align*}
where $\hat{G}$ is the Laplace transform of $ C \mathrm{e}^{A \hat{t}} B$.
\ilcc{The first inequality follows from the Cauchy-Swartz inequality. Note also that}
$\hat{G}$ is strictly proper, \kas{due to the structure of \eqref{sys1}, \eqref{sys2},} with all poles on the open left half plane \kas{(from Assumption \ref{assum_cont_system}(i))}, and hence \ic{its $\mathcal{H}_2$-norm} is finite \kas{(e.g. \cite[Ch. 2]{doyle2013feedback})}.
 Hence, noting that $\epsilon$ can be chosen to be arbitrarily small and that for any $\epsilon_2$, there exists finite $\hat{t}$ such that $|C \kas{\mathrm{e}^{A(t - \tau_1)}} x(\tau_1)| \leq \epsilon_2$
 for all $t \geq \hat{t}$, it follows that for any $\hat{\epsilon} > 0$, there exists $\tilde{\tau}_1 \in \mathbb{R}_+$ such that $|\kas{y(t)}| \leq \epsilon_2 + \delta := \hat{\epsilon}$ for $t \in [\tilde{\tau}_1, \tau_1 + \hat{\tau}_1]$.
 Finally note \ank{that when the value of $\hat{\tau}_1$ is sufficiently large,  it holds that \ka{$\tilde{\tau}_1 + \overline{\tau}_1 \leq \tau_1 + \hat{\tau}_1$}}. \ank{The latter completes the proof.}
 \hfill $\blacksquare$

\emph{Proof of Theorem \ref{thm_conv}:}
\ilc{The trajectories $z(t,j)$ of system \eqref{sys4} are in general non-unique.
  However, it can be trivially shown that for each trajectory of $d^s(t,j)$, there exists a unique trajectory for $\omega(t,j)$, since $\omega$ is the output of linear system \eqref{sys_linear_description} with input $d^s$.}
  The analysis below \ilc{concerns $\omega(t,j)$ given} any \iccc{trajectory $d^s(t,j)$ that is compatible with \eqref{sys4} such that the conditions of Theorem \ref{thm_conv} hold.}
For simplicity, in the analysis below we drop the element $j$ from the argument of the solutions,
\ank{i.e. denoting $x(t,j)$  and $d^c_i(t,j), i \in N$, by simply $x(t)$ and $d^c_i(t), i \in N$ respectively.}


  For system \eqref{sys_linear_description},  from  any initial condition $x(0) \in \mathbb{R}^{n+1}$, $\omega(t)$ is given by
\begin{equation}\label{omega_trajectory}
\omega(t) = C \mathrm{e}^{At} x(0) + \int_0^t C \mathrm{e}^{A\tau} B u(t - \tau) d\tau,
\end{equation}
which suggests that the magnitude of $\omega(t)$ satisfies
\begin{equation}\label{omega_trajectory_magn}
|\omega(t)| \leq |C \mathrm{e}^{At} x(0)| + |\int_0^t C \mathrm{e}^{A\tau} B u(t - \tau) d\tau|.
\end{equation}

Since $A$ is Hurwitz, it follows that for any $\kas{\hat{\epsilon}} > 0$ there exists $\tau \in \mathbb{R}_+$ such that $|C \mathrm{e}^{At} x(0)| \leq \kas{\hat{\epsilon}}$ for all $t~\geq~\tau$.
Furthermore, for the integral part of \eqref{omega_trajectory_magn}, it holds that $|\int_0^t C \mathrm{e}^{A\tau} B u(t~-~\tau) d\tau| \leq  \int_0^t |C \mathrm{e}^{A\tau} B| d\tau \norm{u}_{\infty} \leq
\int_0^\infty |C \mathrm{e}^{A\tau} B| d\tau \norm{u}_{\infty} = \hat{L}\norm{u}_{\infty}$,
\kas{noting that $\hat{L}$ is bounded from Assumption \ref{assum_cont_system}(i).}
 Hence, for any $x(0) \in \mathbb{R}^{n+1}$ and any $\epsilon > 0$, there exists $\tau \in \mathbb{R}_+$ such that $|\omega(t)| \leq \epsilon + \hat{L}\norm{d^s - d^{s,*}}_{\infty} =  \hat{\omega}$ for all $t \geq \tau$.

Now for given $\hat{\omega}$ consider the sets $S(\hat{\omega}) = \{j : \omega^1_j \leq \hat{\omega} \}$ and
  $\hat{S}(\hat{\omega}) = \iccc{N \setminus S(\hat{\omega})}$, which should be interpreted as the sets of loads with and without active frequency feedback. In particular, \ank{since $|\omega(t)| \leq \hat{\omega}, t \geq \tau$, the dynamics of $d^c_j$ reduce from \eqref{sys_temperature}, \eqref{sys_hysteresis_3},} to $\eqref{sys_hysteresis}, \eqref{sys_temperature},$ for $j \in \hat{S}(\hat{\omega})$.
   Furthermore, note that Corollary \ref{corollary_variance} applies to the set $\hat{S}(\hat{\omega})$, suggesting that $\lim_{|N| \rightarrow \infty} \mathbb{V}(d^s_{\hat{S}(\hat{\omega})}) = 0$.

\ank{In the arguments below the variables $\tau_1, \hat{\tau}_1$ and $\tilde{\tau}_1, \overline{\tau}_1$ are used as in Lemmas \ref{prop_dc_sum_convergence} and \ref{lemma_g2} respectively.}
 From Lemma \ref{prop_dc_sum_convergence}, it follows that as $|N| \rightarrow \infty$,  then for any $\kas{{\epsilon}_1} > 0, \hat{\tau}_1 \in \mathbb{R}_+$  there exists $\tau_1 \in \mathbb{R}_+$ such that $\int_{\tau_1}^{\tau_1 + \hat{\tau}_1} (d^s_{\hat{S}(\hat{\omega})}(t) - d^{s,*}_{\hat{S}(\hat{\omega})})^2 dt \leq \kas{{\epsilon_1}}$.
 \kas{Note that the value of $\tau_1$ depends on  $\hat{\tau}_1, \epsilon_1$ and the initial conditions.}
 It then follows by applying Lemma \ref{lemma_g2} with $u = (d^s_{\hat{S}(\hat{\omega})}(t) - d^{s,*}_{\hat{S}(\hat{\omega})})$,
that for any $\kas{\tilde{\epsilon}} > 0, \overline{\tau}_1 \in \mathbb{R}_+$, there exists $\tilde{\tau}_1 \in \mathbb{R}_+$ such that
\kas{$|\omega(t)| \leq
|C \mathrm{e}^{At} x(0)| + |\int_0^t C \mathrm{e}^{A\tau} B u(t - \tau) d\tau|
= |C \mathrm{e}^{At} x(0)|$ $+ |\int_0^t C \mathrm{e}^{A\tau} B (d^s_{\hat{S}(\hat{\omega})}(t-\tau) - d^{s,*}_{\hat{S}(\hat{\omega})}) d\tau| + |\int_0^t C \mathrm{e}^{A\tau} B (d^s_{{S}(\hat{\omega})}(t-\tau) - d^{s,*}_{{S}(\hat{\omega})}) d\tau|
\leq \tilde{\epsilon} + \hat{L}\|d^s_{{S}(\hat{\omega})} - d^{s,*}_{{S}(\hat{\omega})} \|_{\infty}$ for all $t \in [\tilde{\tau}_1, \tilde{\tau}_1 + \overline{\tau}_1]$.}
\kas{Note that, as follows from the arguments in the proof of Lemma \ref{lemma_g2}, it holds that $\tilde{\tau}_1 \geq \tau_1$.}
The rest of the proof is split in two parts, depending on \icll{whether} ${S}(\hat{\omega}) = \emptyset$ or not.

\emph{Part 1:}
If ${S}(\hat{\omega}) = \emptyset$ then the proof is complete from the above arguments.

\emph{Part 2:}
If ${S}(\hat{\omega}) \neq \emptyset$, then from Design condition \ref{des_condition_freq_thresholds} it holds that $\|d^s_{{S}(\hat{\omega})} - d^{s,*}_{{S}(\hat{\omega})}\|_{\infty} \leq  \sum_{j \in S(\hat{\omega})} \zeta_j \overline{d}_j  \leq \max (\hat{L}^{-1}(\hat{\omega} - \delta), 0)$.
Then, letting $\kas{\bar{\epsilon}} \in (0, \delta)$, it follows that $|\omega(t)| \leq  \hat{\omega} - (\delta - \kas{\bar{\epsilon}}) = \hat{\omega}_1$  for all $t \in [\tilde{\tau}_1, \tilde{\tau}_1 + \overline{\tau}_1]$.
Then, note that when $|\omega(t)| \leq \hat{\omega}_1$ the set of loads with active frequency feedback reduces to $S(\hat{\omega}_1)$ which satisfies $|S(\hat{\omega}_1)| \leq |S(\hat{\omega})|$.


The rest of the proof repeats the above argument to construct a decreasing sequence of $\hat{\omega}_i$, \kas{where the subscript $i$ corresponds to the $i$th element of the sequence}.
 In particular, since $|\omega(t)| \leq  \hat{\omega}_1$  for all $t \in [\tilde{\tau}_1, \tilde{\tau}_1 + \overline{\tau}_1]$ it holds that for any $\epsilon_2 > 0, \hat{\tau}_2 \in \mathbb{R}_+$  there exists $\tau_2 \in \mathbb{R}_+$ such that $\int_{\tau_2}^{\tau_2 + \hat{\tau}_2} (d^s_{\hat{S}(\hat{\omega}_1)}(t) - d^{s,*}_{\hat{S}(\hat{\omega}_1)})^2 dt \leq \epsilon_2$ and hence for any $\overline{\tau}_2 \in \mathbb{R}_+$ there exists $\tilde{\tau}_2 \in \mathbb{R}_+$ such that $|\omega(t)| \leq  \hat{\omega}_1 - (\delta - \epsilon) = \hat{\omega}_2$  for all $t \in [\tilde{\tau}_2, \tilde{\tau}_2 + \overline{\tau}_2]$.
 The latter follows from  Lemma \ref{prop_dc_sum_convergence} and Lemma \ref{lemma_g2} as above.

\kas{Below, we define $\tau_i, \hat{\tau}_i$ and $\tilde{\tau}_i, \overline{\tau}_i$ in analogy to $\tau_1, \hat{\tau}_1$ and $\tilde{\tau}_1, \overline{\tau}_1$ corresponding to the $i$th iteration of the considered sequence.
It then follows that the values of $\hat{\tau}_i$ and $\overline{\tau}_i$ can be selected at each iteration such that $[\tilde{\tau}_i, \tilde{\tau}_i + \overline{\tau}_i] \subseteq [\tau_i, \tau_i + \hat{\tau}_i]$ and $\tilde{\tau}_{i+1} \geq \tilde{\tau}_{i}$.}
 Therefore, given that there exists $\hat{\omega}_i$ and $\tilde{\tau}_i$ such that $|\omega(t)| \leq \hat{\omega}_i, t \in [\tilde{\tau}_i, \tilde{\tau}_i + \overline{\tau}_i]$, where $\overline{\tau}_i$ can be arbitrarily large, then there exist $\hat{\omega}_{i+1}$ and $\tilde{\tau}_{i+1} > \tilde{\tau}_{i}$ such that $|\omega(t)| \leq \hat{\omega}_i - (\delta - \epsilon) = \hat{\omega}_{i+1}$ for \iccc{all $t \in [\tilde{\tau}_{i+1}, \tilde{\tau}_{i+1} + \overline{\tau}_{i+1}]$.} 
 \ankk{Note also  \ilc{that $\bar \tau_i, \bar \tau_{i+1}$}
 can be appropriately selected such that $[\tilde{\tau}_{i+1}, \tilde{\tau}_{i+1} + \overline{\tau}_{i+1}] \subseteq [\tilde{\tau}_{i}, \tilde{\tau}_{i} + \overline{\tau}_{i}]$. }

Hence, there exists a decreasing sequence of $\hat{\omega}_i$ such that $0 \leq \hat{\omega}_{i+1} \leq \hat{\omega}_i - (\delta - \kas{\bar{\epsilon}})$ and  $|S(\hat{\omega}_{i+1})| \leq |S(\hat{\omega}_i)|$. Furthermore, there exists some finite $n$ such that $\hat{\omega}_n < \omega_m(N)$ which implies that $|S(\hat{\omega}_n)| = \emptyset$.
 Then, Lemma \ref{prop_dc_sum_convergence} holds for the set $N$ and hence it follows that for any $\epsilon > 0$, there exists $\tilde{\tau}_n$ such that  the trajectories of $\omega$ satisfy $|\omega(t)| \leq \epsilon$ for $t \in [\tilde{\tau}_n, \tilde{\tau}_n + \overline{\tau}_n]$, where $\overline{\tau}_n$
can be selected to be arbitrarily~large.
%
 \hfill $\blacksquare$

\balance

\bibliography{andreas_bib}

\end{document}